\documentclass[11pt,reqno]{article}

\title {Mixing Times of the biased card shuffling and the asymmetric exclusion process}

\author{Itai Benjamini \and Noam Berger \and Christopher Hoffman
\and Elchanan Mossel}

\usepackage{latexsym}
\usepackage{amsfonts}
\usepackage{amsmath}
\usepackage{amsthm}
\usepackage{amssymb}

\newtheorem {theorem}{Theorem}[section]
\newtheorem{claim}[theorem]{Claim}
\newtheorem{thm}[theorem]{Theorem}
\newtheorem{defn}[theorem]{Definition}
\newtheorem{lemma}[theorem]{Lemma}
\newtheorem{cor}[theorem]{Corollary}

\newtheorem{prob}[theorem]{Problem}

\newcommand{\E}{\mbox{$\bf E$}}

\newcommand{\C}{C}
\newcommand{\CC}{C}
\newcommand{\D}{D}
\newcommand{\K}{K}
\newcommand{\Z}{{\mathbb Z}}

\newcommand{\CA}{{\mathcal{CA}}}
\newcommand{\EXZ}{{\mathcal{EX}}(\Z,p)}
\newcommand{\EX}{{\mathcal{EX}}}
\newcommand{\EXS}{{\mathcal{EX}_2}(\Z,p)}

\newcommand{\var}{\mbox{$\bf Var$}}

\renewcommand{\prob}{{\bf P}}
\renewcommand{\P}{{\bf P}}
\newcommand{\ve}{\text{ and }}
\newcommand{\ixigul}{{0 \hspace{-.06in}{\sf X}}}

\def\eps{\epsilon}
\def\lam{\lambda}

\begin{document}

\maketitle

\begin{abstract}
Consider the following method of card shuffling. Start with a deck
of $N$ cards numbered 1 through $N$. Fix a parameter $p$ between 0
and 1. In this model a ``shuffle'' consists of uniformly selecting
a pair of adjacent cards and then flipping a coin that is heads
with probability $p$.  If the coin comes up heads then we arrange
the two cards so that the lower numbered card comes before the
higher numbered card.  If the coin comes up tails then we arrange
the cards with the higher numbered card first. In this paper we
prove that for all $p \neq 1/2$, the mixing time of this card
shuffling is $O(N^2)$, as conjectured by Diaconis and Ram
\cite{DR}. 
Our result is a rare case of an exact estimate for the
convergence rate of the Metropolis algorithm.
A novel feature of our proof is that the analysis of an
infinite (asymmetric exclusion) process plays an essential role in
bounding the mixing time of a finite process.

\end{abstract}

\section{Introduction}

The Metropolis algorithm is a widely used algorithm for sampling 
distributions on large finite sets. A variety of techniques were
developed in order to analyze the convergence rate of the Metropolis
algorithm (see \cite{DS}), 
yet in many problems arising in applications we do not know how to
estimate the convergence rate. In this paper we introduce new
techniques which allow the
analysis of a card shuffling Metropolis algorithm which was non-amenable to
the standard techniques in the field. 

Card shuffling procedures provide a natural family of Markov
chains which played a crucial role in the development of the
theory of the convergence rate of Markov chains, see e.g.
\cite{AD,Dia88,Dia98}. In this paper we analyze the mixing time of
a {\em biased} card shuffling procedure which has a nonuniform
stationary distribution.


\subsection{The biased card shuffling chains}
As the biased card shuffling has a finite state space, we
formulate it as discrete time Markov chain.
\begin{defn} \label{def:ca}
For $0 \leq p \leq 1$, let {\bf $\CA_d(N,p)$} denote the following
discrete time Markov chain on permutations of $N$ cards labeled
$1,\ldots,N$. A step of the chain consists of selecting uniformly
at random  a pair of adjacent cards and then flipping a coin that
is heads with probability $p$. If the coin comes up heads then we
arrange the two cards so that the lower numbered card comes before
the higher numbered card. If the coin comes up tails then we
arrange the cards with the higher numbered card first.
\end{defn}

\noindent
Note that if $p=1/2$ then the stationary distribution for
$\CA_d(N,p)$ is the uniform distribution on $S_N$
(the set of all permutations on $N$ elements),
but if $p \neq 1/2$ then the stationary distribution is not uniform.

A novel feature of our results is that the heart of the proof of the mixing result for
the finite card shuffling model involves the analysis of infinite processes on $\Z$.
Since infinite processes are
naturally defined in continuous time, we use a
continuous time version of the biased card shuffling model.

\begin{defn}
For $0 \leq p \leq 1$, let {\bf $\CA(N,p)$} denote the following
continuous time Markov chain on permutations of $N$ cards labeled
$1,\ldots,N$. Each pair of adjacent cards $i,i+1$ is picked with
rate $1$ independently. Then we toss a coin which is heads with
probability $p$. If the coin comes up heads then we arrange the
two cards so that the lower numbered card comes before the higher
numbered card. If the coin comes up tails then we arrange the
cards with the higher numbered card first.
\end{defn}

Diaconis and Ram \cite{DR} were interested in the following
slightly different chain:
\begin{defn}
For $0.5 \leq p \leq 1$, let $q=1-p$ and let $\theta=q/p$. The
{\bf Metropolis biased card shuffling} is the following discrete
time Markov chain on permutations of $N$ cards labeled
$1,\ldots,N$. A step of the chain starts with selecting uniformly
at random a pair of adjacent cards. If the two cards are arranged
in a decreasing order, then we switch them. If they are arranged
in an increasing order, then with probability $\theta$ we switch
them and with probability $1-\theta$ we do nothing.
\end{defn}

\subsection{Main results}

 The {\bf total-variation} distance between measures $\mu$ and
$\nu$ on a finite space $X$ is
$$||\mu-\nu||_{TV}=\frac12 \sum_{x \in X}
|\mu(x)-\nu(x)|=\sup_{A \subset X}|\mu(A)-\nu(A)|.$$
 Now we define
the mixing time of a Markov chain $\sigma$ on a finite state space
$X$.  For any $x \in X$ let $x_t$ be the distribution on $X$ at
time $t$ under the action of $\sigma$. The {\bf mixing time} of
the Markov chain is defined by
$$\tau_1= \inf\{t:\ \sup_{x,x' \in X} ||x_t -x'_t||_{TV} \leq e^{-1}\}.$$

Our main result is
\begin{thm} \label{thm:disccards}
For all $p \neq 1/2$ there exists a constant ${\K} = {\K}(p)$ such
that the mixing time of the discrete time biased card shuffling on
$N$ cards is at most ${\K} N^2$.
\end{thm}

\begin{cor}
For all $p > 1/2$ there exists a constant ${\K} = {\K}(p)$ such
that the mixing time of the Metropolis biased card shuffling on
$N$ cards is at most ${\K} N^2$.
\end{cor}
\begin{proof}
The discrete time card shuffling is a slow down of the Metropolis
card shuffling. More precisely, consider the following process: At
every time, with probability $1-p$ we do nothing, and with
probability $p$ we do a Metropolis step. This process is the
discrete time card shuffling.
\end{proof}

\noindent
This verifies a conjecture of Diaconis and Ram \cite{DR}. A lower
bound for the mixing time of the form $N^2$ is easy and well
known.

As the only difference between the two processes is that the
continuous time process is ``$N-1$ times faster'' the following
result is equivalent.
\begin{thm} \label{thm:cards}
For all $p \neq 1/2$ there exists a constant ${\K} = {\K}(p)$ such
that the mixing time of the continuous time biased card shuffling
on $N$ cards is at most ${\K} N$.
\end{thm}
As our proofs are all done with continuous time processes, we will
prove Theorem \ref{thm:cards} and derive Theorem
\ref{thm:disccards} as a corollary.

\subsection{Motivations and related results}
When running the Metropolis algorithm for sampling distribution on
large finite sets, it is necessary to insure that the algorithm
converges rapidly. Various techniques were developed in order to bound
the convergence rate (see Subsection \ref{subsec:remarks}), yet in
many cases none of these methods apply. 
Such an example is the asymmetric card shuffling chain (see Subsection
\ref{subsec:remarks}) which we analyze in a novel way in this paper. 

Our result also suggests an interesting comparison between the 
``systematic scan'' and the ``random scan'' heuristics in sampling
(see e.g. \cite{Fishman} or \cite{DR}).

In \cite{DR} Diaconis
and Ram studied a different version of biased card shuffling.  In
their model the selection of the pair of adjacent cards was not
random, but done in a prescribed deterministic manner (``systematic
scan''). 
Their discrete time model, like ours (``random scan''), 
has a mixing time of $O(N^2)$ for $p \neq 1/2$. 
Our result may be interpreted as saying that the ``systematic scan''
doesn't give an improvement over the ``random scan''.

In \cite{FPRU} the authors introduce a model of computation where each comparison
operation has probability $p < 1/2$ of returning the true result and probability $1-p$
of returning a false result independently of other comparisons.
The chain ${\mathcal CA}_d(N,1-p)$ (${\mathcal CA}_d(N,p)$) is
performing the randomized version of bubble sort in this noisy computation model.
Our result shows the robustness to noise of the randomized bubble sort algorithm as the convergence
time of ${\mathcal CA}_d(N,1-p)$ is $O(N^2)$ for all $p > 1/2$.

We would also like to remark that asymmetric exclusion
processes, which are the key tool in our proof, also play a
crucial role in the study of the quantum Heisenberg model.

We conclude this subsection by discussing some of the history of card
shuffling problems.
Gilbert, Shannon and Reed began the mathematical study of card
shuffling by introducing a good model for how people shuffle cards
\cite{GS,R}.
The celebrated theorem of Bayer and Diaconis \cite{BD} states that for the
Gilbert-Shannon-Reed model of card shuffling it takes seven
shuffles in order for a standard 52 card deck to be well mixed.
More generally, \cite{BD} proved that for an
$N$ card deck the mixing time for the Gilbert-Shannon-Reed model
is approximately $\frac32 \log_2 N$.

In the wake of Bayer and Diaconis's result there have been a
number of articles analyzing the mixing time for various methods
of card shuffling. Most relevant to this paper are results of
Wilson as well as Diaconis and Ram.
Wilson \cite{wilson} found that the mixing
time for ${\mathcal CA}_d(N,1/2)$ to within a factor of 2.
The mixing time is of order $N^3 \log N$.
Note the sharp contrast with Theorem
\ref{thm:disccards} where we show if $p\neq 1/2$ then the mixing
time of ${\mathcal CA}_d(N,p)$ is $O(N^2)$. 

\subsection{The asymmetric exclusion process}
Most of the proof of our main result is devoted to analysis of
the asymmetric exclusion processes. We now define these processes
which are of independent interest, see e.g. \cite{liget1,liget2}.
First we define our family of finite exclusion processes.  The
process $\EX(N,k,p)$ will be an exclusion process with $N$
containers and $k$ particles.

\begin{defn} \label{def:ex}
Let $k$ and $N$ be integers such that $1 \leq k < N$, and $0 \leq
p \leq 1$. Let $\EX(N,k,p)$ be the continuous time Markov process
defined on
$$X_{n,k}=\left \{x \in \{0,1\}^{[1,N]}:\ \sum_{i=1}^N x_i = k \right \}$$
in the following way. Given the current state $x$, each pair of
coordinates $i,i+1$ of $x$ is picked at rate $1$. If
$x_i = x_{i+1}$, then the chain will stay at state $x$; otherwise,
the two coordinates $i,i+1$ will be reassigned as $(x_i,x_{i+1}) =
(1,0)$ with probability $p$, and as $(x_i, x_{i+1}) = (0,1)$ with
probability $1-p$.
\end{defn}

\begin{defn} \label{def:exZ}
Let $0 \leq p \leq 1$. Let $\EX(\Z,p)$ be the continuous time
Markov process defined on $\{0,1\}^{\Z}$ in the following way.
Given the current state $x$, each pair of coordinates $i,i+1$ of
$x$ is picked at rate $1$. If $x_i = x_{i+1}$, then the
chain will stay at state $x$; otherwise, the two coordinates
$i,i+1$ will be reassigned as $(x_i,x_{i+1}) = (1,0)$ with
probability $p$, and as $(x_i, x_{i+1}) = (0,1)$ with probability
$1-p$.
\end{defn}

We are particularly interested in the set
$$A=\left \{a:\ \sum_{-\infty}^{-1}(1-a_i)=\sum_{0}^{\infty}a_i<\infty.\right \}$$
There is a partial order on the set. We write $a \succeq b$ if for
all $r$
\begin{equation} \label{eq:partial order}
\sum_{i=-\infty}^r{1-a_i} \leq \sum_{i=-\infty}^r{1-b_i}.
\end{equation}

The maximal state in $A$ is the ground state
\begin{equation} \label{eq:mZr}
G_{\Z}(i)= \left\{
        \begin{array}{ll}
                1               &i<0\\
                0               &i\geq 0\\
        \end{array}
\right.
\end{equation}

The aspect of asymmetric exclusion processes that we are most
interested in is the tail of the  hitting time. Given any $x \in
A$ (or measure $\mu$ on $A$) the {\bf hitting time, $H(x)$ (or
$H(\mu)$),} is defined by
\begin{equation} \label{eq:H}
H(x)=\inf\{t:\ x_t=G_{\Z}\}
\end{equation}

In particular we want to consider $H(I_N)$ where
\begin{equation} \label{eq:INr}
I_{N}(i)= \left\{
        \begin{array}{ll}
                1               &i< -N\\
                0               &i\in [-N,-1]\\
                1               &i\in [0,N-1]\\
                0               &i\geq N\\
        \end{array}
\right.
\end{equation}

\begin{thm}\label{oneovern}
For all $p > 1/2$ and $\eps > 0$ there exists a constant
${\D}={\D}(p,\eps)$ s.t.
\begin{equation*}
\prob(H(I_N)<DN)>1-\frac{\eps}{N}.
\end{equation*}
\end{thm}
In Section \ref{sec:height} we show that Theorem \ref{oneovern}
implies Theorem \ref{thm:cards}. Most of the work in this paper is
in proving Theorem \ref{oneovern}.

\subsection{Remarks on analytic techniques} \label{subsec:remarks}
When applying standard analytic techniques to the study of the
biased card shuffling problem we encounter several problems which
prevent us from obtaining sharp bounds for the mixing time. We discuss briefly the
difficulties in estimating the mixing time in terms of the
spectral gap and the log Sobolev constant. The results of
\cite{DR} also suggest that a group theoretic approach is hard to
apply for solving this particular problem.

The standard bound for the mixing time in terms of the spectral
gap of the generator of the Markov process, $ -\lam_2$, (see
\cite{saloff} for background) yields
\begin{equation} \label{eq:spectral_gap}
\tau \leq \frac{2}{-\lam_2} \log \frac{1}{\min_x \pi(x)},
\end{equation}
where $\pi$ is the stationary distribution.
Moreover, combining our reduction from the card shuffling
to the exclusion process with the bound on the spectral gap of the (continuous time)
exclusion process in \cite{CM}, it is straightforward to verify
that there exists positive $c_1$ and $c_2$ such that indeed the spectral
gap satisfies $c_1 \leq -\lam_2(N) \leq c_2 $ for all $N$.
However, the probability
space contains elements of very small probability, so the term
$\log (1 / \min_x \pi(x))$ is of order $N^2$ (see \cite{DR} where
the stationary distribution for the Metropolis chain is given).
Thus (\ref{eq:spectral_gap}) yields a bound of order $N^2$.

A standard way to reduce the dependency on the smallest
probability, is to use the log Sobolev constant $\alpha$ instead
of the spectral gap with the estimate
\begin{equation} \label{eq:log_sob}
\tau \leq  \frac{4}{\alpha} \log_{+} \log \frac{1}{\min_x \pi(x)}.
\end{equation}
However, plugging the indicator of the set which consists of a
single element $(N \ldots 1)$ in the variational formula of the
log Sobolev constant (see \cite{saloff}) implies that (for the
continuous time model) $\alpha = O(1/N^2)$. (We use the notation
$f(N)=O(N)$ $(f(N)=\Omega(N))$ if there exists a constant
$c<\infty$ ($c>0$) such that $f(N)<cN$ ($f(N)>cN$).) Thus
(\ref{eq:log_sob}) doesn't give the right bound of $N$ on the
mixing time.

\subsection{Road Map}
We conclude the introduction with an overview of the main steps of
the proof.
\begin{enumerate}
\item
In Section \ref{sec:height} we show how Theorem \ref{oneovern}
implies Theorems \ref{thm:cards} and \ref{thm:disccards}. The
reduction follows \cite{wilson} in using {\bf height functions},
together with coupling arguments. The height functions provide a
coupling of the biased card shuffling to the finite exclusion
processes.  This reduces the problem of bounding the mixing time
for the biased card shuffling to bounding the tail of the hitting
times of some finite exclusion process. Then we couple the finite
exclusion processes with the infinite exclusion process. This
reduces the problem to bounding the tail of $H(I_N)$. Now we can
use some of the machinery developed for the study of exclusion
processes on $\Z$.

\item
In Section \ref{sec:motionless} we define the blocking measure
$\Psi$ on $\{0,1\}^{\Z}$.  We show that $\Psi$ is an invariant
measure for $\EX(\Z,p)$. In Section
\ref{sec:main_idea} we will see how to bound the tail of $H(I_N)$
in terms of the tail of $H(\Psi)$.

\item
In Section \ref{sec:second} we introduce an asymmetric exclusion
process with second class particles.  Second class particles are a
common tool used in exclusion processes (see. e.g. \cite{liget2}).
This will be our main tool for proving Theorem \ref{oneovern}. We
discuss some of the processes related to the asymmetric exclusion
process with second class particles.

\item
In Section \ref{sec:main_idea} we use exclusion processes with
second class particles to bound the tail of $H(I_N)$ in terms of
the tail of $H(\Psi)$. Then we
bound the tail of $H(\Psi)$. This bound
allows us to prove Theorem \ref{oneovern}, which in turn allows us
to bound the mixing time of the biased card shuffling process.

\end{enumerate}

\section{Coupling card shuffling to exclusion processes} \label{sec:height}

In this section we show how Theorem \ref{oneovern} implies Theorem
\ref{thm:cards}.  This reduces bounding the mixing time for the
biased card shuffling to bounding the tail of $H(I_N)$.
Following \cite{wilson}, we use the following collection of height
functions to map a permutation of a deck of cards to an exclusion
process configuration.

For any $k,\ 1\leq k<N,$ consider the map $h_k: S_N \to X_{N,k}$
defined by
\begin{equation*}
(h_k(\pi))_i= \left\{
        \begin{array}{ll}
                1               &\pi_{i} \leq k\\
                0               &\pi_{i} > k
        \end{array}
\right.
\end{equation*}
It is easy to see that
\begin{claim} \label{claim:det}
$\pi$ is determined by $(h_k(\pi))_{k=1}^{N-1}$.
\end{claim}

For $\pi \in
S_N$, we write {\bf $\pi_t$} for the random
variable representing the value of the process at time  $t$ that
starts at $\pi$ and evolves according to $\CA(N,p)$. Similarly for
$x \in X_{N,k}$  (or $\{0,1\}^{\Z}$) we let the random variable $x_t$
represent the configuration at time $t$ for the exclusion process
that started at $x$ and evolves according to $\EX(N,k,p)$
($\EXZ$).

\begin{claim}
For all $N$ and $k$, the processes $h_k(\pi_{t})$ are Markovian.
Moreover, $h_k(\pi_{t})$ evolves according to  the exclusion
process $\EX(N,k,p)$.
\end{claim}

Throughout the paper we will define a number of couplings.  The
most important of these we refer to as {\bf canonical couplings}.
The idea of the canonical couplings is that we have a collection
of initial conditions of permutations (or exclusion process
states) and we use one set of clocks and one set of biased coin
flips to update all of the process simultaneously.

We begin by defining a coupling of the process $\CA(N,p)$ for all the
configurations $\pi \in S_N$ in the following way: A transition
of shuffling is to be performed by choosing a pair of adjacent
coordinates $i,i+1$ at rate $1$, and tossing a coin $X$ which is
heads with probability $p$. If $X=H$ then we rearrange the cards
in coordinates $i,i+1$ in increasing order, while if $X=T$ then we 
rearrange the cards in coordinates $i,i+1$ in decreasing order.
The same pair of coordinates $i,i+1$ and the same coin $X$ are chosen for all
$\pi \in S_n$ simultaneously.
We call this coupling the {\bf canonical coupling} for $\CA(N,p)$.

We can similarly define a coupling for all the configurations in $\EX(N,k,p)$.
A transition is to be performed by choosing a pair of adjacent coordinates $i,i+1$ in $\Z$ at rate $1$ and tossing a
coin $X$ which is heads with probability $p$. For a state $x$ of $\EX(\Z,p)$ we will update $x$ as follows.
If $x_i = x_{i+1}$, then we do nothing. Otherwise, if $X=H$, we let $x_i = 0, x_{i+1} = 1$, and if $X=T$, we let
$x_{i+1} = 0, x_i = 1$. Again, the same pair of coordinates and the same coin
is chosen for all states of $\EX(N,k,p)$. We call this coupling the 
{\bf canonical coupling} for $\EX(N,k,p)$. In the same way we define a 
canonical coupling for $\EX(\Z,p)$.

It is immediate to verify that
\begin{claim} \label{claim:canonical}
For all $N,k$ and $p$, the canonical couplings for $\CA(N,P),\EX(N,k,p)$ and
$\EX(\Z,p)$ are all well defined and have the right marginals.

Moreover for all $N,k$ and $p$, the map $h_k$ maps the canonical coupling
of $\CA(N,p)$ to the canonical coupling of $\EX(N,k,p)$, i.e., if
$(\pi_t)_{\pi \in S_N, t \geq 0}$ evolves according to the canonical 
coupling for $S_N$, then the process 
\[
\left( \{h_k(\pi_t) : \pi \in S_N\} \right)_{t \geq 0}
\]
has the same distribution as the process 
\[
\left( \{\sigma_t : \sigma \in X_{N,k}\} \right)_{t \geq 0}
\]
where $(\sigma_t)_{\sigma \in X_{N,k}, t \geq 0}$ evolve according to the 
canonical coupling for $\EX(N,k,p)$.
\end{claim}

The analysis of the process $\EX(N,k,p)$ utilizes monotonicity properties
of the canonical coupling. For $a,b \in X_{N,k}$ we write $a \succeq b$, if
for all $r$, $\sum_{i=1}^r a_i \geq \sum_{i=1}^r b_i$.
The maximal state with respect to this partial order is
\[
g_{N,k}(i)= \left\{
        \begin{array}{ll}
                1               &i \leq k\\
                0               &i > k\\
        \end{array}
\right. ,
\]
and the minimal state with respect to this partial order is
\[
m_{N,k}(i) = \left\{
        \begin{array}{ll}
                0               &i \leq k\\
                1               &i > k\\
        \end{array}
\right. .
\]
We let $H(N,k)$ be the hitting time of 
the state $g_{N,k}$ for the process $\EX(N,k,p)$ started at $m_{N,k}$.


It is immediate to see that
\begin{claim} \label{claim:ex}
The canonical couplings for $\EX(N,k,p)$ and $\EX(\Z,p)$ 
are monotone.  That is, for both processes if $x \succeq y$,
then for all $t$ it holds that $x_t \succeq y_t$.
\end{claim}

Since $g_{N,k}$ and $m_{N,k}$ are the maximal and minimal elements with
respect to $\succeq$ it follows that
\begin{claim} \label{claim:ebound}
Under the canonical coupling for $\EX(N,k,p)$ it holds that
\[
\prob(\exists x,y \in X_{N,k} \mbox{ s.t. } x_t \neq y_t) 
= \prob(H(N,k) >t).
\]
\end{claim}

\begin{lemma} \label{shuffling bound}
Under the canonical coupling for $\CA(N,p)$ it holds that
\[
\prob(\exists \sigma,\tau \in S_N \mbox{ s.t. } \sigma_t \neq \tau_t) 
\leq \sum_{k=1}^{N-1} \prob(H(N,k) > t).
\]
\end{lemma}

\begin{proof}
\begin{eqnarray} \label{eq:det}
\prob(\exists \sigma,\tau \in S_N \mbox{ s.t. } \sigma_t \neq \tau_t) &=& 
\prob(\exists 1 \leq k \leq N-1, \sigma,\tau \in S_N \mbox{ s.t. } 
h_k(\sigma_t) \neq h_k(\tau_t)) \\ 
&\leq& 
\sum_{k=1}^{N-1}  \prob(\exists x,y \in X_{N,k} \mbox{ s.t. } x_t \neq y_t) 
= \sum_{k=1}^{N-1} \P[H(N,k) > t], \label{eq:ebound}
\end{eqnarray}
where (\ref{eq:det}) follows from Claim \ref{claim:det} and (\ref{eq:ebound})
follows from claim \ref{claim:canonical} and \ref{claim:ebound}.
\end{proof}

The remaining coupling step is to couple the finite processes to an 
infinite process. 
\begin{lemma} \label{finite_infinite}
For all $p \geq 1/2$, all $N$, all $1 \leq k < N$, 
and all $t > 0$ the processes
$\EX(N,k,p)$ and $\EX(\Z,p)$ satisfy that
\[
\prob(H(N,k) > t) \leq \prob(H(I_N) > t).
\] 
\end{lemma}

\begin{proof}
Consider the map $X_{N,k} \to \{0,1\}^{\Z}$ 
sending $x \in X_{N,k}$ to $\widehat{x}$ with
\begin{equation}\label{tilde}
\widehat{x}(i)=\left\{
        \begin{array}{ll}
                1           &   i<-k\\
                x(i+k+1)    &   i \in [-k,N-k-1]\\
                0           &   i \geq N-k
        \end{array}
\right.
\end{equation}

We will now couple $\EX(N,k,p)$ to the process $\EX(\Z,p)$.
More formally, we will couple the canonical coupling of
$\EX(N,k,p)$ with the canonical coupling of $\EX(\Z,p)$.

For the process $\EX(\Z,p)$ we pick a pair of coordinates $i,i+1$ at
rate $1$ and then use a coin $X$ to rearrange the coordinates $i,i+1$ in the
usual manner. 
Now, if $1 \leq i+k+1 \leq i+k+2 \leq N$, then we use the same coin
$X$ to rearrange the coordinates $i+k+1,i+k+2$ for the process $\EX(N,k,p)$.

Clearly this coupling is well defined and has the right marginals.
It is easy to see that if $p \geq 1/2$, then this coupling has the following
important property: For all $x \in X_{n,k}$ and all
$t \geq 0$, it holds that $\widehat{(x_t)} \succeq (\widehat{x})_t$.

Writing $I = I_N, g = g_{N,k}$, and $m = m_{N,k}$, and since 
$\widehat{m} \succeq I$, it follows 
that under this coupling, if $I_t = G_{\Z}$, then  
\[
\widehat{(m_t)} \succeq (\widehat{m})_t \succeq I_t = G_{\Z}.
\]
and therefore $m_t = g$. The claim of the lemma follows. 
\end{proof}

%
%

We end this section by noting that the canonical coupling of the
biased card shuffling shows that with high probability by time
$DN$ all of the processes agree.

\begin{lemma} \label{shuffling bound2}
Theorem \ref{oneovern} implies that the canonical coupling for $\CA(N,p)$
has
\begin{equation} \label{eq:canon}
\P(\exists \sigma, \tau \in S_N \mbox{ s.t. } \sigma_{D N} \neq \tau_{D N}) 
< \eps.
\end{equation}
In particular Theorem \ref{oneovern} 
implies theorems \ref{thm:cards} and \ref{thm:disccards}.
\end{lemma}

\begin{proof}
By Theorem \ref{oneovern}, $\P[H(I_N) > D N] < \eps/N$, and therefore by
Lemma \ref{finite_infinite} for all $k$ it holds that 
$\P[H(N,k) > D N] < \eps/N$. It now follows from 
Lemma \ref{shuffling bound} that 
\[
\P(\exists \sigma, \tau \in S_N \mbox{ s.t. } \sigma_{D N} \neq \tau_{D N}) < 
(N-1) \eps / N < \eps,
\]   
so we obtain (\ref{eq:canon}). Taking $\eps = e^{-1}$, we deduce
Theorem \ref{thm:cards} which immediately implies Theorem 
\ref{thm:disccards} 
\end{proof}

\section{The blocking measure}
\label{sec:motionless}

In this section we define a distribution $\Psi$ on $\{0,1\}^{\Z}$
which is invariant under the action of $\EX(\Z,p)$.
$\Psi$ is known as the {\bf blocking measure} (see e.g.
\cite{liget1}). In Section \ref{sec:main_idea} we will bound the
tail of $H(\Psi)$ and show that the tail of $H(I_N)$ can be
bounded in terms of the tail of $H(\Psi)$.

Fix any $p>1/2$. Define $\mu=\mu(p)$ on $\{0,1\}^{\Z}$, to be the
product measure with probabilities
\begin{equation} \label{eq:skew_prod}
\mu(\eta(i)=1)= \left(\frac{1-p}{p}\right)^i
\left/\left(1+\left(\frac{1-p}{p}\right)^i \right)\right.
\end{equation}

\begin{lemma} \label{block} \cite{liget1}
The measure $\mu$ is stationary for $\EXZ$.
\end{lemma}
\begin{proof}
This is proven on page 381 of \cite{liget1}.
\end{proof}
Notice that $\mu$ is supported on configurations $\eta$ s.t. there
exists a $C_\eta$ s.t. $\eta(i)=1$ for every $i<-C_\eta$ and
$\eta(i)=0$ for every $i>C_\eta$. There are only countably many
configurations of this type and each of them has a positive
measure. We have already defined
$$A=\left \{a:\ \sum_{-\infty}^{-1}(1-a_i)=\sum_{0}^{\infty}a_i<\infty.\right \}$$

\begin{defn} \label{def:motionless}
The blocking measure {\bf $\Psi$} on $\{0,1\}^{\Z}$ is defined by
\begin{equation}\label{legal}
\Psi=\mu|_{A}.
\end{equation}
\end{defn}

\begin{cor}
$\Psi$ is stationary and ergodic for the exclusion process.
\end{cor}

By Poincare's recurrence theorem and the fact that
$\Psi(G_{\Z})>0$, we get the following lemma.
\begin{lemma}\label{recurr}
$\lim_{T \to \infty}\prob(H(\Psi)>T)=0$
\end{lemma}
In Lemma \ref{tochelet} we will show that
$\prob(H(\Psi)>T)=e^{-\Omega(\sqrt{T})}.$

\begin{lemma} \label{psitail}
\begin{equation}\label{zanav_shel_psi}
\Psi(\exists_{i>N}(\eta(i)=1))=\Psi(\exists_{i<-N}(\eta(i)=0))
=O\left( \left( \frac{1-p}{p}\right)^N \right).
\end{equation}
For any $T>0$
\begin{equation}\label{zanav_shel_psi2}
\prob(\exists t \in (T,T+N) \mbox{ and $i>2N$ such that }
    \eta_t(i)=1)=e^{-\Omega(N)}
\end{equation}
\end{lemma}

\begin{proof}
In the product measure, the $\mu$ probability that there exists an
occupied site right of position $N$ is bounded by
\begin{eqnarray*}
\sum_{i=N+1}^{\infty} \frac{((1-p)/p)^i}{1+((1-p)/p)^i}
\leq\sum_{i=N+1}^{\infty}\left(\frac{1-p}{p}\right)^i =O\left(
\left( \frac{1-p}{p}\right)^N \right).
\end{eqnarray*}
Since $\Psi$ is obtained from the product measure by conditioning
on an event of positive probability the first part of the lemma is
true.  For the second part if there exists $t \in (T,T+N)$ and
$i>2N$ such that $\eta_t(i)=1$ then either
\begin{enumerate}
\item there exists $i\geq N$ such that $\eta_T(i)=1$ or
\item for some $t$
    $$\max\{i:\eta_t(i)=1\}-\max\{i:\eta_T(i)=1\}>N$$
\end{enumerate}
By the first part of the lemma the probability of the first event
is decreasing exponentially in $N$.  The second event happens only
if the right most particle moves to the right $N$ times in a
period of length $N$.  As moves to the right happen with rate
$1-p<\frac12$ the probability that this happens is also decreasing
exponentially in $N$.
\end{proof}

\section{Exclusion processes with second class particles}
\label{sec:second}

The main tool that we will use in the rest of the paper is adding
second class particles to our exclusion process.  We now describe
some of the basics about exclusion processes with second class
particles.  For a more rigorous treatment of exclusion processes
with second class particles see
\cite{liget2}.

\begin{defn} \label{def:exS}
Let $0 \leq p \leq 1$. Let $\EXS$ be the continuous time Markov
process defined on $\{0,1,2\}^{\Z}$ in the following way. Given
the current state $x$, each pair of coordinates $i,i+1$ of $x$ is
picked at rate $1$. If $x_i = x_{i+1}$, then the chain will stay
at state $x$. If the two coordinates $i,i+1$ initially are $(0,1)$
or $(1,0)$ then they will be reassigned as $(x_i,x_{i+1}) = (1,0)$
with probability $p$, and as $(x_i, x_{i+1}) = (0,1)$ with
probability $1-p$. If initially they are $(0,2)$ or $(2,0)$ then
they will be reassigned as $(2,0)$ with probability $p$, and as
$(0,2)$ with probability $1-p$. If initially they are $(1,2)$ or
$(2,1)$ then they will be reassigned as $(1,2)$ with probability
$p$, and as $(2,1)$ with probability $1-p$.
\end{defn}
If $x_i=1$ then we say that there is a first class particle in
position $i$, if $x_i=2$ then we say that there is a second
class particle in position $i$, and if $x_i = 0$, then we say that
the site $i$ is empty.

It is helpful to have in mind the following ordering of $0$, $1$ and $2$.
Particle $1$ has priority over $0$ and $2$ in moving the left. Particle $2$ has
priority over $0$ (but not over $1$) in moving to the left. Particle of type $2$ is therefore
ranked in-between particle of type $0$ and particle of type $1$.

It is therefore natural to consider to the following two
projections. In the first projection $\delta^{2 \to 1}$, $2$'s are
projected to $1$'s, while in the second projection $\delta^{2 \to
0}$, $2$'s are projected to $0$'s. More formally,

\begin{equation*}
\delta_t^{2 \to 1}(i)= \left\{
        \begin{array}{ll}
                0               &\delta_t(i)=0, \\
                1               &\delta_t(i)>0.
        \end{array}
\right.
\end{equation*}
and
\begin{equation*}
\delta_t^{2 \to 0}(i)= \left\{
        \begin{array}{ll}
                0               &\delta_t(i)\neq 1, \\
                1               &\delta_t(i)=1.
        \end{array}
\right.
\end{equation*}

\begin{claim} \label{claim:evolve}
Both $\delta^{2 \to 1}$ and $\delta^{2 \to 0}$ evolve according to $\EXZ$.
\end{claim}

The next process we consider represents the dynamics between
particles of type $1$ and particles of type $2$ and eliminates all
the information on the $0$'s.

To define $\delta^{0 \hspace{-.06in}{\sf X}}_t$ we first eliminate all of the zeroes from
$\delta_t$, and then change all of the twos to zeroes. This is
only well defined up to translation, so we must also decide which
translate we want.  We do this by tagging one particle in
$\delta_t$ and having $\delta^{0 \hspace{-.06in}{\sf X}}_t(0)$ correspond to the tagged
particle.

We now make this more formal. Let
$$u_0(0) = \left\{
    \begin{array}{ll}
        \sup\{i:\ \delta_0(i)=1\} & \mbox{if } \sup\{i:\delta_0(i)=1\}<\infty,\\
        \sup\{i< 0:\ \delta^1_0(i)=1\} & \mbox{otherwise.}
    \end{array}
\right.$$
We refer to the particle which is in position $u_0(0)$ at time 0
as the {\bf tagged particle}. Let $u_t(0)$ be the location of the
tagged particle at time $t$. For $n=1,2,...$, let
$$u_t(n)=\min\{i>u_t(n-1):\delta_t(i)>0\},$$ and
$$u_t(-n)=\max\{i<u_t(-n+1):\delta_t(i)>0\}.$$
Thus $u_t(n)$ represents the location of the $n$th particle at
time $t$. Let

\begin{equation}\label{defdelta4}
\delta^{0 \hspace{-.06in}{\sf X}}_t(i)=\left\{
        \begin{array}{ll}
                0               &\delta_t(u_t(i))=2\\
                1               &\delta_t(u_t(i))=1
        \end{array}
\right.
\end{equation}

\begin{lemma} \label{lem:stochdom}
If $\delta$ evolves according to $\EXS$ and has initial
distribution $\delta^{0 \hspace{-.06in}{\sf X}}_0$ stochastically
dominating $\Psi$, then $\delta^{0 \hspace{-.06in}{\sf X}}_t$
stochastically dominates $\Psi$ for all $t$.
\end{lemma}

\begin{proof}
In \cite{liget1} it is shown that for all $i \in \Z$, $\Psi$ is invariant under the Markov operator on
$\{0,1\}^{\Z}$ which at rate $1$ tosses a coin $X$ which heads with probability $p$, then if $x_i \neq  x_{i+1}$ and $X=H$
then $x_i$ and $x_{i+1}$ are updated as $x_i = 1, x_{i+1} = 0$, if $x_i \neq x_{i+1}$ and $X=T$, then $x_i$ and $x_{i+1}$
are updated as $x_i = 0, x_{i+1} = 1$. This implies that $\Psi$ is an invariant measure for the process
$\delta^{0 \hspace{-.06in}{\sf X}}$.
It now immediately follows that if $\delta^{0 \hspace{-.06in}{\sf X}}_0$ stochastically dominates $\Psi$, then
$\delta^{0 \hspace{-.06in}{\sf X}}_t$ stochastically dominates $\Psi$ for all $t$.
\end{proof}




\section{Proof of Main Results}\label{sec:main_idea}

Let $\{Y_i\}_{i\in\Z}$ be i.i.d. random variables s.t.
$\prob(Y_i=0)=\prob(Y_i=1)=1/2$, and let $Z_i=2Y_i$. The main tool
to prove Theorem \ref{oneovern} is to study $\EXS$ with initial
conditions
\begin{equation} \label{eq:stoch_INr5}
\sigma_0(i)= \left\{
        \begin{array}{ll}
                1               &i< -N\\
                0               &i\in [-N,-1]\\
                1               &i\in [0,N-1]\\
                Z_i             &i\geq N
        \end{array}
\right. ,
\end{equation}
This is useful for proving Theorem \ref{oneovern} because
$\sigma^{2\to 0}_t =(I_N)_t$.

For any $a \in \{0,1\}^{\Z}$ such that
 $\lim_{i \to -\infty}a_i=1$ we set
\begin{equation}\label{eq:L}
L(a)=\min\{i:\ a_i=0\}.
\end{equation}
This indicates the left most empty position. The same way, for $a$
s.t. $\lim_{i \to \infty}a_i=0$, we indicate the right most
particle by
\begin{equation*}
R(a)=\max\{i:\ a_i=1\}.
\end{equation*}

For a constant ${\C}$ we define three events:
\begin{equation} \label{A1}
A_1(C,N)=\{\forall t \in ({\C} N, ({\C}+1)N) \ L(\sigma^{2\to
1}_t)>2N \},
\end{equation}
\begin{equation} \label{A2}
A_2(C,N)=\{ \forall t \in ({\C} N,({\C}+1)N) \
R(\sigma^\ixigul_t)<2N \},
\end{equation} and
\begin{equation} \label{A3}
A_3(C,N)=\{\exists t \in ({\C} N, ({\C}+1)N)
    \mbox{ such that }\sigma^\ixigul_t=G_{\Z}\}.
\end{equation}

\begin{lemma} \label{lem:prob}
\begin{equation*}
\prob  \left(H(I_N)\leq ({\C}+1)N\right)
\end{equation*}
\begin{equation*}
\geq 1- \prob \left(A_3^c(C,N)|A_1(C,N),A_2(C,N)\right) -
\prob\left(A_1^c(C,N)\right)-\prob\left(A_2^c(C,N)\right).
\end{equation*}
\end{lemma}

\begin{proof}
Recall that by $G_{\Z}$ we denote the ground state (equation
(\ref{eq:mZr})). Notice that $\sigma^{2\to 0}_t=G_{\Z}$ if
$\sigma^\ixigul_t=G_{\Z}$ and $L(\sigma^{2\to 1}_t)>0$.  Thus if
$A_1(C,N)$ and $A_3(C,N)$ both occur then there exists
$t\leq(\C+1)N$ such that $\sigma^{2\to 0}_t=(I_N)_t=G_{\Z}$.
\end{proof}

Although the previous lemma did not depend on the definition of
$A_2(C,N)$, we use it because it is easy to bound
 $$\prob \left(A_3^c(C,N)|A_1(C,N),A_2(C,N)\right)$$ in terms of the tail of $H(\Psi)$.

\begin{lemma} \label{condprob}
$$\prob \left(A_3^c(C,N)|A_1(C,N),A_2(C,N)\right)\leq \prob
(H(\Psi)>N).$$
\end{lemma}

\begin{proof}
If $A_1(C,N)$ and $A_2(C,N)$ both happen then $\sigma^\ixigul$
behaves according to $\EXZ$ conditioned on the event that there is
never a particle to the right of $2N$. By Lemma \ref{lem:stochdom}
the distribution of $\sigma^\ixigul_{{\C}N}$ stochastically
dominates $\Psi$. Putting these two facts together gives us
 $\prob\left(A_3^c(C,N)|A_1(C,N),A_2(C,N)\right)\leq \prob(H(\Psi)>N)$.
\end{proof}

It is also not difficult to bound $\prob\left(A_2^c(C,N)\right)$.

\begin{lemma} \label{lem:A2}
$\prob\left(A_2^c(C,N)\right)=e^{-\Omega(N)}$.
\end{lemma}

\begin{proof}
This follows from Lemmas \ref{psitail} and \ref{lem:stochdom}.
\end{proof}

Our next goal is to bound $\prob(A_1^c(C,N))$. Then we will bound
the tail of $H(\Psi)$.

In order to bound $\prob(A_1^c(C,N))$ we first bound $\prob(\tilde
A_1^c(C,N))$ where
$$\tilde A_1(C,N)=\{L(\sigma^{2\to 1}_{{\C}N})>3N\}$$
and use the following lemma.

\begin{lemma} \label{lem:Atilde}
$\prob(A_1^c(C,N)) \leq \prob(\tilde A_1^c(C,N)) +e^{-\Omega(N)}$.
\end{lemma}

\begin{proof}
If $\tilde A_1^c(C,N)$ happens but $A_1^c(C,N)$ does not, then the
leftmost container with out a particle moves to the left at least
$N$ times in time $N$. For that to happen, the clock left of this
container has to ring $N$ times. But since its rate is smaller
than $1$, by simple large deviation estimates for Poisson
variables, the probability of this happening is decreasing
exponentially in $N$.
\end{proof}

To bound the probability of $\tilde A_1^c(C,N)$  we study the
process $\beta$ which has initial distribution
\begin{equation} \label{eq:stoch_INr2}
\beta_0(i)= \left\{
        \begin{array}{ll}
                Y_i             &i \leq 0\\
                Z_i=2Y_i               &i > 0
        \end{array}
\right. .
\end{equation}

I.e. the initial distribution of $\beta$ is the following: Every
place contains a particle with probability $1/2$, and the places
are independent of each other. Left of the origin the particles
are first class particles, while right of the origin the particles
are second class particles.

We are interested in the processes $\beta^{2\to 0}$ and
$\beta^{2\to 1}$. The process $\beta^{2\to 1}$ is the stationary
i.i.d. process. The process $\beta^{2\to 0}$, on the other hand,
is the process that starts with no particles on the right half of
the line , and an i.i.d. measure on its left half.

Let $x(t)$ be the location of the tagged particle in $\beta^{2\to
0}$ at time $t$ and let $x'(t)$ be the location of the tagged
particle in $\beta^{2\to 1}$ at time $t$. We will bound the
expectation and variance of $x(t)$.



The following lemma follows from the proof of the shock wave
phenomenon in \cite{liget2}. For the convenience of the reader, we
prove it here again.

\begin{lemma}\label{couple_with_hetzi}

There exists $\varrho<1$ such that under the canonical coupling,
for every $n$ and for every time $t$,
\begin{equation*}
\prob(|x(t)-x'(t)|>n)<\varrho^n.
\end{equation*}
\end{lemma}
\begin{proof}
We define the processes $\beta^\ell$ and $\beta^d$ ($\ell$ stands
for locations and $d$ stands for distances):

$\beta^\ell_t(i)$ is the location of the $i$-th particle in
$\beta^{2\to 1}$. To be more precise, $\beta^\ell_t(0)=x'(t)$,
and, inductively, for positive $i$ we take
$$\beta^\ell_t(i)=\min\left(j:j>\beta^\ell_t(i-1)\ve\beta^{2\to 1}_t(j)=1
\right)$$ and for negative $i$, equivalently, we take
$$\beta^\ell_t(i)=\max\left(j:j<\beta^\ell_t(i+1)\ve\beta^{2\to 1}_t(j)=1
\right).$$ $\beta^d_t(i)$ is defined to be
$\beta^\ell_t(i)-\beta^\ell_t(i-1)$.

Since for every $t$, $\left\{\beta^{2\to 1}_t(i)\right\}_{i\in\Z}$
is distributed according to the $(1/2,1/2)$ product measure, we
get that for every $t$, $\left\{\beta^d_t(i)\right\}_{i\in\Z}$ are
i.i.d. geometric variables with parameter $1/2$.

Let
$$s(t)=\sup\{i:\beta^\ixigul(i)=1\}.$$
Then for all $t$
$$x(t)-x'(t)=\sum_{i=1}^{s(t)}\beta^d_t(i).$$
By Lemmas \ref{psitail} and \ref{lem:stochdom}
$$\P(s(t)> n/3 )
 =\Psi\left( \exists i>  n/3: \eta(i)=1 \right)
 =O\left(\frac{1-p}{p}\right)^{n/3}.$$
 By the distribution of $\left\{\beta^d_t(i)\right\}_{i\in\Z}$ there
exists $\alpha<1$ such that
$$\P\left(\sum_1^{ n/3  }\beta^d_t(i)>n\right)=O(\alpha^n).$$
Thus there exists $\varrho<1$ such that
$$ \prob(|x(t)-x'(t)|>n)
 \leq \P(s(t)> n/3 )+ \P\left(\sum_1^{ n/3  }\beta^d_t(i)>n\right)
 <\varrho^n $$
\end{proof}

The following lemma is proved in \cite{kipnis}:
\begin{lemma}[Kipnis]\label{kipnis}
\begin{eqnarray*}
\lim_{t\to\infty}\frac{\E(x'(t))}{t}=s'=-\frac12(2p-1).
\end{eqnarray*}
\begin{eqnarray*}
\lim_{t\to\infty}\frac{\var(x'(t))}{t}=v' \leq \frac{6}{2p-1}
<\infty.
\end{eqnarray*}
\end{lemma}

\noindent {\bf Remark}: note that $s<0$.

Combining Lemma \ref{kipnis} and Lemma \ref{couple_with_hetzi} and
the fact that $\E(x(t))$ and $\var(x(t))$ are continuous in $t$,
we get
\begin{lemma}\label{shock_loc}
There exists $t_0$ and there exist $s'<s<0$ and $v'<v<\infty$ s.t.
for every $t\geq t_0$,
\begin{equation*}
\frac{\E(x(t))}{t}<s.
\end{equation*}
and
\begin{equation*}
\frac{\var(x(t))}{t}<v.
\end{equation*}
\end{lemma}

Consider the exclusion process $\gamma$ which has the initial
distribution

\begin{equation} \label{eq:stoch_INr4}
\gamma_0(i)= \left\{
        \begin{array}{ll}
               1                   &i < 0\\
               1-Y_{i}            &i \geq 0
        \end{array}
\right. .
\end{equation}

One can couple a copy $\gamma'$ of $\gamma$ with $\beta^{2\to 0}$
such that for all $t$ and $i$
$$\gamma'_t(i)=1-\beta^{2\to 0}_t(-i).$$
$x(t)$ corresponds to the location of the rightmost particle in
$\beta^{2\to 0}_t$, so the processes $-x(t)$ and $L(\gamma_t)$
have the same law.


Applying Chebychev's inequality we get the following estimate:

\begin{lemma} \label{cheb}
For any $\delta>0$ and any $t>t_0$
$$\prob\left(x(t)\leq\left(s+\sqrt{\frac{v}{\delta}}\right)t\right)>1-\delta/t$$
and
$$\prob\left(L(\gamma_t)\geq-\left(s+\sqrt{\frac{v}{\delta}}\right)t\right)>1-\delta/t.$$
\end{lemma}

For any $l$ we define $\gamma^l$ to be the process starting at
\begin{equation*}
\gamma^l_0(i)= \left\{
        \begin{array}{ll}
               1                   &i < -l\\
               1-Y_{i}            &i \geq -l
        \end{array}
\right. .
\end{equation*}

We get from Lemma \ref{cheb} that for any $t>t_0$,
\begin{equation}\label{gammal}
\prob\left(L(\gamma^l_t)\geq-l-\left(s+\sqrt{\frac{v}{\delta}}\right)t\right)
=\prob\left(L(\gamma_t)\geq-\left(s+\sqrt{\frac{v}{\delta}}\right)t\right)>1-\delta/t.
\end{equation}

Now we are ready to bound $\prob({\tilde A}_1^c)$.

\begin{lemma} \label{dhf}
For any $\eps > 0$, there exists a constant ${\C} = {\C}(p,\eps)$
such that for all $N\geq[t_0]+1$
\begin{equation*}
\prob(\tilde A_1^c(C,N))<\frac{\eps}{N}.
\end{equation*}
\end{lemma}

\begin{proof}
As the canonical coupling preserves domination, if
$$\gamma^l_0 \succeq \sigma^{2\to 0}_0 \mbox{ and } L(\gamma^l_t)>3N$$ then
$L(\sigma^{2\to 0}_t)>3N.$ This gives us that for any $l$ and $t$
\begin{equation} \label{lCdelta}
  \prob(L(\sigma^{2\to 0}_{t})>3N)
\geq 1-\prob(\gamma^l_0 \not\succeq \sigma^{2\to 0}_0)
        -\prob(L(\gamma^l_t)\leq 3N)
\end{equation}

Choose $j$ such that for all $N$
$$\prob\left(\sum_{i=1}^{jN}Y_i<N\right)<\epsilon/2N.$$
Then the lemma follows from equation (\ref{gammal}) and equation
(\ref{lCdelta}) with
    $l,\delta,{\C}$ and $t$ chosen such that
        $l=jN$,
        $\delta=4v/s^2$,
        ${\C}>\max(-2(3+j)/s,2\delta/\eps,1)$, and
        $t={\C}N$.

This is because
$$\prob\left(\gamma^{jN}_0 \not\succeq \sigma^{2\to 0}_0\right)
    =\prob\left(\sum_{i=1}^{jN}Y_i<N\right)<\eps/2N$$
and by Lemma \ref{cheb} and equation \ref{gammal},
\begin{eqnarray*}
 \prob(L(\gamma^{jN}_{{\C}N})> 3N)
   & =& \prob(L(\gamma_{{\C}N})> (3+j)N)\\
   & >& \prob(L(\gamma_{{\C}N})>-\frac{s}{2}{\C}N)\\
   & >& \prob(L(\gamma_{{\C}N})>-\left(s+\sqrt{\frac{v}{\delta}}\right){\C}N)\\
   & >& \delta/ {\C}N \\
   & >& 1-\eps/2N.
\end{eqnarray*}

\end{proof}

\begin{lemma}\label{almost}
For every $N>[t_0]+1$ and $\eps > 0$,
\begin{equation}\label{eq:almost}
\prob(H(I_N)<(C+1)N)
    >1- \eps/N-\prob(H(\Psi)>N).
\end{equation}
\end{lemma}

\begin{proof}
This follows from Lemmas \ref{lem:prob}, \ref{condprob},
\ref{lem:A2}, \ref{lem:Atilde},  and \ref{dhf}.

\end{proof}

In order to prove Theorem \ref{oneovern} we first prove the following
lemma:
\begin{lemma}\label{tochelet}
\begin{equation*}
\prob(H(\Psi) \geq N)=e^{-\Omega(\sqrt{N})}
\end{equation*}
\end{lemma}
\begin{proof}
For every $N$ large enough, we wish to estimate the probability
that $H(\Psi)\geq ({\C}+1)N^2$ where ${\C}$ is the constant from
Lemma \ref{almost}. We take $N$ large enough so that the
probability in (\ref{eq:almost}) is bigger than $\frac{1}{2}$.
Such $N$ exists by Lemma \ref{recurr}. Recall that
\begin{equation*}
I_N(i)= \left\{
        \begin{array}{ll}
                1               &i < -N\\
                0               &i\in [-N,-1]\\
                1               &i\in [0,N-1]\\
                0               &i\geq N
        \end{array}
\right.
\end{equation*}
Now, for every $j=0,1,2,...,N$, let $P_j=\prob(H(\Psi)\geq
(C+1)Ni)$. Of course, $P_0=1$. Now, we proceed inductively. Let
$$U_N=1-\prob(\eta_t \succeq I_N).$$
Notice that by Lemma \ref{lem:stochdom} it does not depend on $t$,
and, by (\ref{zanav_shel_psi}), $$U_N=e^{-\Omega(N)}.$$ For every
$t$,
$$\prob(H(\Psi)\geq t+(C+1)N| \eta_t \succeq I_N)<\frac{1}{2}.$$
Therefore,
$P_i\leq\frac{P_{i-1}}{2}+U_N$ for every $i>0$. Therefore,
\begin{equation*}
\prob(H(\Psi)\geq (C+1)N^2)=P_N \leq 2^{-N} + \sum_{i=1}^N 2^{-i}
U_N = e^{-\Omega(N)}.
\end{equation*}
By monotonicity we can interpolate and get that for every $t$
\begin{equation*}
\prob(H(\Psi)\geq (C+1)t)=e^{-\Omega(\sqrt{t})}
\end{equation*}
and thus
\begin{equation*}
\prob(H(\Psi)\geq t)=e^{-\Omega(\sqrt{t})}.
\end{equation*}
\end{proof}

We can now prove our main results.
\begin{proof}[Proof of Theorems \ref{thm:disccards}, \ref{thm:cards}, and \ref{oneovern}]
By Lemma \ref{tochelet},
$H(\Psi)=e^{-\Omega(\sqrt{N})}=o(N^{-1})$. Therefore, by
(\ref{eq:almost}),
\begin{equation*}
\prob(H(I_N)<(C(p,\eps)+1)N)>1-\frac{\eps}{N}-o(N^{-1}).
\end{equation*}
Taking ${\D}=\C(p,\frac{\eps}{2})+1$ we get Theorem \ref{oneovern}
is satisfied for all arbitrarily large $N$.  Thus we can choose
${\D}$ so that it is true for all $N$. Theorems
\ref{thm:disccards} and \ref{thm:cards} follow by Lemma
\ref{shuffling bound}.
\end{proof}

We conclude the paper with a brief comment about how
${\D}={\D}(p)$ depends on $p$. We see that in Lemma \ref{kipnis}
that $s=-\frac12(2p-1)$ and $v\leq
\frac{6}{2p-1}$. For large $N$, in Lemma \ref{shock_loc} using
$\eps=1/e$ we can choose
 $${\CC}= \frac{8ve}{s^2}\leq\frac{1024e}{(2p-1)^{3}}.$$
For large $N$ we can choose
 $${\D}=2{\CC} \leq\frac{2048e}{(2p-1)^{3}}.$$
It is easy to show that ${\D}$ must be chosen bigger than
$1/(2p-1).$ The discrepancy in the power of $2p-1$ comes from the
use of Chebychev's inequality in Lemma \ref{cheb}.
We believe that a more careful analysis would allow one to choose
$D$ s.t. $D=\theta(1/(2p-1))$.

\medskip\noindent
{\bf Acknowledgment:} We thank Persi Diaconis for important comments on
a previous version of this paper.
We thank Dror Weitz, Prasad Tetali,
Thomas Liggett and David Aldous for many interesting discussions. Part of the
work was done while Noam Berger was an intern in Microsoft
Research.  Christopher Hoffman was partially supported by NSF
grant \#0100445.

\end{document}